\documentclass[notitlepage, 11pt]{article}

\usepackage{color}
\usepackage{amsmath}
\usepackage{latexsym,amssymb}
\newtheorem{lemma}{Lemma}
\newtheorem{theorem}{Theorem}

\newtheorem{definition}{Definition}

\newtheorem{proposition}{Proposition}

\newtheorem{remark}{Remark}

\newcommand{\noi}{\noindent}
\newcommand{\spa}{\vspace{.2in}}

\newcommand{\SINR}{\mathrm{SIR}}
\newcommand{\SINNR}{\mathrm{SINR}}
\newcommand{\Exp}{\mathbb{E}}
\newcommand{\R}{\mathbb{R}}

\newcommand{\calG}{\mathcal{G}}

\newcommand{\calM}{\mathcal{M}}
\newcommand{\Pb}{\mathbb{P}}

\newcommand{\Z}{\mathbb{Z}}
\newcommand{\I}{\mathcal{I}}

\newcommand{\bfp}{\mathbf{p}}

\newcommand{\eps}{\varepsilon}

\usepackage{latexsym}
\usepackage{amsmath, amssymb, amsfonts}

\numberwithin{equation}{section}

\title{On Scaling Limits of Power Law Shot-noise Fields
\footnote{This work was supported by an award from Simons Foundation ({\#}197982) to the University of Texas at Austin.}}

\author{Fran\c{c}ois Baccelli and Anup Biswas\\ [2mm]
Department of Electrical and Computer Engineering,\\
University of Texas, Austin.}

\date{November 19, 2014}

\begin{document}

\maketitle

\begin{abstract}
This article studies the scaling limit of a class of shot-noise
fields defined on an independently marked stationary Poisson
point process and with a power law response function.
Under appropriate conditions, it is shown that the shot-noise
field can be scaled suitably to have a non degenerate $\alpha$-stable limit, 
as the intensity of the underlying point process goes to infinity.
More precisely, finite dimensional distributions are shown to converge 
and the finite dimensional distributions
of the limiting random field have i.i.d. stable random components.
We hence propose to call this limit the $\alpha$- stable white noise field.
Analogous results are also obtained for the extremal shot-noise field
which converges to a Fr\'{e}chet white noise field.
Finally, these results are applied to the modeling and
analysis of interference fields in large wireless networks. 
\vspace{.2in}


\vspace{.2in}

\noi{\bf Keywords:}\ Poisson point process, shot-noise process, percolation, $\alpha$-stable white noise, scaling limits.

\end{abstract}

\section{Introduction}
The present paper is focused on scaling limits for a class of
shot-noise fields \cite{SKM} associated with stationary
Poisson point processes and with a power law response function.
It is motivated by the modeling of ultra-dense wireless networks and
the analysis of interference fields that arise in this type of
stochastic networks \cite{Baccelli-Bartek}.
By scaling limits, we understand the analogue of a functional
central limit theorem, namely a rescaling of the field such that,
in the limit, this field has non-degenerate joint distributions.
To give a first example, let $L$ be a bounded, non-negative,
radial function (i.e. $L(x)=L(\|x\|)$) on $\R^d$ that is integrable
 w.r.t. to Lebesgue measure $m(dy)$ on $\R^d$.
Now if $\Phi_\lambda$ is a stationary Poisson point process on
$\R^d$ of intensity $\lambda$, the corresponding shot-noise field is given by
$$\I_\lambda(L;z)=\sum_{\{x_i\}\in\Phi_\lambda}L(\|x_i-z\|).$$
This field is translation invariant.
It turns out that $\I_\lambda(L;0)$ has finite moments.
In fact, $\Exp[\I_\lambda(L;0)]=\lambda\int_{\R^d}L(\|y\|)\, m(dy)$.
Therefore we can define the second order approximation of
$\I_\lambda$ as follows
$$\hat{\I}_\lambda(L;z)=\frac{1}{\sqrt{\lambda}}\Big(\I_\lambda(L;z)-\Exp[\I_\lambda(L;0)]\Big).$$
Then using the Laplace transform, it is not difficult to
show that as $\lambda\to\infty$, the scaled field
$\hat{\I}_\lambda$ converges to a Gaussian random field.
In particular, given two points $z_1, z_2\in \R^d$,
$(\hat{\I}_\lambda(L;z_1), \hat{\I}_\lambda(L;z_2))$ converges
to a $2$-dimensional Gaussian random variable with
covariance matrix given by 
$$\begin{pmatrix}
\int_{\R^d}L^2(\|z_1-y\|)\, m(dy) & \int_{\R^d}L(\|z_1-y\|)L(\|z_2-y\|)\, m(dy)\\

\int_{\R^d}L(\|z_1-y\|)L(\|z_2-y\|)\, m(dy) & \int_{\R^d}L^2(\|z_2-y\|)\, m(dy)
\end{pmatrix}.$$

It turns out that the response function that
is the most commonly used in the wireless literature is
a power law \cite{Baccelli-Bartek} i.e.,
$L(x)=\frac{1}{\|x\|^\beta}, \beta>d$, so that the central limit
scale is not the appropriate one. Two interesting limits come
from the identification of the appropriate scaling: 
the {\em stable white noise field} and the
{\em Fr\'echet white noise field}.
A particularly nice property of these fields is their
independence property: all their finite dimensional
distributions are those of random vectors with independent
and identically distributed components. In general,
this property is not shared by fields for which the right
scale is the central limit one (see e.g. \cite{bierme, breton, kaj}).
The convergence of the random fields are
defined in terms of the convergence of their
finite dimensional distributions. 
In the present paper, the underlying point process is a
Poisson point process (PPP) and we use the Laplace transform
to establish convergence. Scaling limits of random fields
fed by a PPP have already been thoroughly studied in the literature.
In the $1$-dimensional setting, scaling properties of Poisson
shot-noise processes were thoroughly studied in
\cite{klup-kuhn-04, klup-mik-03, klup-mik-95}.
For other $1$-dimensional point processes than Poisson,
see also \cite{doney-obrien, hsing-teugels}. In higher dimension,
\cite{hein-molcha} studies limits of union random fields
associated with a stationary point process and
defined by certain classes of functions that are regularly varying at zero.
Properties of extremes of shot-noise fields defined by a 
bounded slowly varying function are studied in \cite{lebedev}.
\cite{dombry} establishes various qualitative properties
of extremal random fields. Let us mention that the main difference
whith these papers is that the class
of functions that are used in the present paper to analyze
high density shot-noise (or extremal shot-noise) fields, are singular
at $0$ and not integrable on $\R^d$. It is shown below that this singularity
at $0$ is instrumental to get the independence (white noise) properties of
the limiting random fields, which are one of the main 
findings of the present paper.

These limit results are discussed
in Section~\ref{secmr}.
Section~\ref{secapp} discusses a few applications
of these limit results to the modeling of communication rates
in ultra-dense wireless networks where interference is treated as noise.
For the detailed context, see \cite{Baccelli-Bartek}.
The networks considered feature wireless transmitters located
according to some realization of a homogeneous PPP
of intensity $\lambda$ in $\R^2$. The interference field
is modeled as the shot-noise field of this Poisson point
process for the response function $r^{-\beta}$. Each transmitter 
maintains a wireless link (an information theoretic channel) to
its receiver, which is assumed to be located at distance 1 from it, in some
random and uniform direction. It is shown in Section~\ref{secapp} that 
the results of Section~\ref{secmr} can be used to obtain
the speed of decrease of the SIR (Signal to Interference Ratio)
of a typical link in such a network when $\lambda$ tends to infinity.
The interest in this question stems from the following result
of information theory: when treating interference as noise,
the (Shannon) communication rate obtained by this typical link 
is proportional to $\log (1+\SINR)$ \cite{CoT91}.
Hence, in the first place, these limiting results allow
one to predict the speed of decrease of the Shannon rate of the
typical link when the network is densified (i.e. $\lambda$ tends
to infinity).

These limiting results also allow one to estimate
the speed at which the SIR decreases with $\lambda$ in a network with density 1
when all links now have length $\lambda$. More precisely, it is shown in 
Section~\ref{secapp} that in such a scenario, it is possible to
transmit $1$ bit over distance $\lambda$ with a delay $D_\lambda$ such that
the ration $\frac{D_\lambda}{\lambda^{\beta}}$ has a non-degenerate limit
with high probability when $d$ tends to infinity (Theorem~\ref{theo-2}).
To compare this result with the existing ones,
we recall a result from \cite{Baccelli-Bartek-Omid}
which can be expressed as follows: for all fixed positive
SIR threshold\footnote{Here it is assumed that
there is threshold on the SIR of a link above (resp. below) which 
communication is assumed successful (resp. unsuccessful).},
for all possibly single or multi-hop strategies,
the expected \textit{delay}, say $\Exp[D_\lambda]$,
to transmit $1$ bit over distance $\lambda$ in a Poisson
wireless network of the type described above and with density $1$
grows faster than $\lambda$, i.e.,
$\frac{\Exp[D_\lambda]}{\lambda}\to\infty$ as $\lambda\to\infty$.
%Finally, in Theorem~\ref{theo-3} we show that
%there is a decay rate of the joint SIR
%(SINR) thresholds for multiple links in a ultra-dense network such that the links behave independently as the intensity grows to infinity.
 
Finally, we use
the joint scaling limit results of Section~\ref{secmr} to 
derive percolation properties of the 
SINR graph \cite{Fran-Mees,Baccelli-Bartek}.
The Signal to Interference and Noise Ratio (SINR)
differs from the SIR defined above as some constant or
random term, called
the thermal noise power,
is added to the interference in the denominator of the ratio of SINR
\cite{Baccelli-Bartek}.
The SINR graph is a random geometric graph with an edge between
two nodes if the SINR from one to the other is above some threshold.
The results of Section~\ref{secmr} are used to estimate 
the speed of decrease to 0 of 
the SINR threshold that makes this graph
percolate when the network 
density tends to infinity (Theorems~\ref{theo-3} and \ref{theo-4}).

\section{Model and Scaling Limits}
\label{secmr}
As mentioned above, one of the most common response functions used in wireless networking
(for $d=2$) is $\frac{1}{\|x\|^\beta}, \beta>2,$ and this is neither bounded
nor globally integrable.
Nevertheless, the additive and the extremal shot-noise fields
generated by this type of response function are finite a.s. 
We are interested in their scaling limits.
%We show that the limiting random field is completely different. In fact we show that if $f$ satisfies certain conditions then any finite dimensional distribution of the limiting field comes from i.i.d. $2/\beta$-stable distributions. Therefore the limiting
%field enjoys an independence property and has analogies with the White Noise field. Similar behavior is also seen for the extremal fields
%generated by $f$. 

Let $m(dy)$ denote the Lebesgue measure on $\R^d$. By $\|\cdot\|$ we denote the Euclidean norm on $\R^d$. We start by defining the Poisson point process of intensity $\lambda$ on $\R^d$.
\begin{definition}
A Poisson point process (PPP) $\Phi_\lambda$ of intensity $\lambda$ on $\R^d$ is a stationary point process such that for any bounded disjoint collection of Borel 
sets $A_i, i=1,\ldots, k,$ of $\R^d$ we have
$$ \Pb\big(\Phi(A_1)=n_1, \Phi(A_2)=n_2,\ldots, \Phi(A_k)=n_k\big)=\Pi_{i=1}^k\Big[e^{-\lambda m(A_i)}\frac{\big(\lambda m(A_i)\big)^{n_i}}{n_i!}\Big].$$
\end{definition}
By independently marked PPP we mean a marked point
process $\tilde\Phi=\{(x_i, \bfp_i)\}_i$ where the locations $\Phi=\{x_i\}$ are given by a PPP and the markings $\{\bfp_i\}$ are independent of the
PPP, i.e., $\Pb(\bfp\in A|\Phi)=\int_A F(dp)$, $A\subset (0, \infty)$ Borel set, is independent of $\Phi$. Marks are assumed to take values in $(0, \infty)$. We also assume that $\Exp[\bfp]=\int p\, F(dp)<\infty$.
By the intensity of a independently marked PPP $\tilde\Phi$ we refer to the intensity
of the underlying PPP $\Phi$.

Now we recall the following classical result on the Laplace functional of an independently marked Poisson point process (see e.g. \cite{Baccelli-Bartek, daley-vere}) which will play a central role in our analysis.

\begin{proposition}\label{laplace}
Let $\tilde\Phi_\lambda$ be an independently marked PPP of intensity $\lambda$ and 
$\mathcal{L}_{\tilde\Phi_\lambda}$ be its Laplace functional.
Then for any non-negative measurable $g$ we have
$$\mathcal{L}_{\tilde\Phi_\lambda}(g)=\exp\Big\{-\lambda\int_{\R^2}\Big(1-\int e^{-g(y, p)} F(dp)\Big) m(dy)\Big\},$$
where $F(\cdot)$ denotes the distribution of the marks.
\end{proposition}
%Note that we have abused the notation of Laplace transformation by denoting it terms of the underlying PPP.
Now we define the \textit{shot-noise} random field corresponding to a marked point process. Given a marked point process $\tilde\Phi=\{(x_i, p_i)\}$ on $\R^d\times (0, \infty)$  and a response function $L:\R^d\times \R^d\times (0, \infty)\to [0, \infty)$, 
the corresponding shot-noise field at $y\in\R^d$ is defined by
$$ \I(L;z)=\sum_{(x_i, \bfp_i)\in\tilde\Phi}L(z, x_i, \bfp_i).$$
Let $f:[0, \infty)\to [0, \infty]$ be a measurable function with
the property that for some positive $\varrho>0$,
\begin{equation}\label{gen-fn}
f(r)=\frac{1}{r^\beta}, \quad \text{for}\ 0\leq r\leq \varrho,\; \; \sup_{r\geq \varrho} f(r)<\infty,
\quad\text{and},\; \int_\varrho^\infty r^{d-1}f(r)dr<\infty,
\end{equation}
for some $\beta>d$.
All the response functions considered in this article
have the form $L(z, x, p)=p\cdot f(\|z-x\|)$ which is commonly
used in the wireless network literature as already explained. We will denote
by $\I(f; \cdot)$ the shot-noise corresponding to the response function $L(z, x, p)=p\cdot f(\|x-z\|)$.
A special response function that will occur often in this article is 
\begin{equation}\label{res-fn}
L(z,x, p)=\frac{p}{\|x-z\|^\beta},
\end{equation}
for some $\beta>d$. We use the notation $\I$ to denote the
shot-noise corresponding to the response function \eqref{res-fn}, i.e., 
$$\I(z)=\sum_{(x_i, \bfp_i)\in\tilde\Phi}\frac{\bfp_i}{\|x_i-z\|^\beta}.$$
Now if $\tilde\Phi$ is an independently marked PPP, then for any $f$
given by \eqref{gen-fn}, the following facts are well known
\cite{Baccelli-Bartek}, Part I, Chapter 2:
\begin{itemize}
\item[-] $\I(f; 0)$ is almost surely finite,
\item[-] $\I(f; \cdot)$ is stationary.
\end{itemize} 
 The finiteness property follows from the integrability
condition of $f$ (last condition of (\ref{gen-fn})).
In fact, if $B$ denotes the unit ball around $0$ then
$$\Exp\Bigl[\sum_{(x_i, \bfp_i)\in\tilde\Phi[B^c]}\bfp_if(\|x_i\|)\Bigr]
= \int_{B^c}f(x)m(dx)\int pF(dp)=
\omega(d)\int_{1}^\infty r^{d-1}f(r) dr\int pF(dp)<\infty,$$
where $\tilde\Phi[A]$ denotes the independently marked PPP with Poisson point
process restricted to the Borel set $A$ and $\omega(d)$ denotes the surface area
of $B$. Therefore finiteness
of $\I(f; 0)$ would follow if one has
\begin{equation}\label{101}
\Pb\Bigl(\sum_{(x_i, \bfp_i)\in\tilde\Phi[B]}\bfp_if(\|x_i\|)<\infty\Bigr)=1.
\end{equation}
Let $\Phi[B]$ be the restriction of PPP to $B$. Then
\begin{align}\label{102}
\Pb\Bigl(\sum_{(x_i, \bfp_i)\in\tilde\Phi[B]}\bfp_if(\|x_i\|)<\infty\Bigr)
&=\sum_{k=1}^\infty
\Pb\Bigl(\sum_{(x_i, \bfp_i)\in\tilde\Phi[B]}\bfp_i f(\|x_i\|)<\infty\big|
\Phi[B]=k\Bigr) \Pb(\Phi[B]=k).
\end{align}
Using the property of Poisson point process for every $k\geq 1$, we can have a i.i.d.
sequence $\{X_1,\ldots, X_k\}$, uniformly distributed on $B$, such that
\begin{equation}\label{103}
\Pb\Bigl(\sum_{(x_i, \bfp_i)\in\tilde\Phi[B]}\bfp_i f(\|x_i\|)<\infty\big|
\Phi[B]=k\Bigr)
= \Pb\Bigl(\sum_{i=1}^k\bfp_i f(\|X_i\|)<\infty\Bigr)=1.
\end{equation}
Hence \eqref{101}  follows from \eqref{102} and \eqref{103}.
Stationarity of $\I(f, \cdot)$ is obvious from the stationary behavior of the underlying Poisson point process.

\subsection{Scaling Limits of Shot-noise Fields}\label{subsec-shot}

We are interested in the limiting behavior of the shot-noise field when the intensity of the underlying PPP goes to infinity. 
By $\I_\lambda(f;\cdot)$ (resp. $\I_\lambda$) we denote the shot-noise field w.r.t.
to an independently marked PPP of intensity $\lambda$ and response function
given by $L(z, x, p)=p\cdot f(\|x-z\|)$ where $f$ satisfies \eqref{gen-fn} (resp (\ref{res-fn})).
Let $\kappa=\frac{\beta}{d}\in(1, \infty)$. Define 
$$\hat{\I}_\lambda(f;\cdot)=\frac{1}{\lambda^\kappa}\I_\lambda(f;\cdot) \quad (\mbox{resp. } 
\hat{\I}_\lambda=\frac{1}{\lambda^\kappa}\I_\lambda).$$
For $n\geq 1$,
and  $(x, t)\in \R^d\times[0, \infty)$ we define 
$$\mathcal{L}_\lambda(f;  x, t)=\Exp\Big[e^{-t\,\hat{\I}_\lambda(f;x)}\Big],\quad
\left(\mbox{resp. } \mathcal{L}_\lambda(x, t)=\Exp\Big[e^{-t\, \hat{\I}_\lambda(x)}\Big]
\right).$$
%Therefore $\mathcal{L}^f_n(t, x)$ denotes the
%Laplace transform for $\hat{\I}^f_n(x)$ at the point $t$.
%To avoid pathological cases, we will assume that $F(\{0\})=0$.
We recall that marks are independently distributed with distribution $F(dp)$ and
$\Exp[\bfp]=\int p\, F(dp)<\infty$.

\begin{lemma}\label{lem1}
Let $x\in\R^d$. Define $\alpha=\frac{1}{\kappa}=\frac{d}{\beta}$.
%If $\mathcal{L}^f_n(x, t)$
%(resp. $\calL_n(x, t)$) denotes the Laplace transform of $\hat{\I}^f_n(x)$ (resp. $\hat{\I}_n(x)$) at the point $t\geq 0$, 
Then, for all $x\in \R^d$ and $t\ge 0$, 
as $\lambda\to\infty$, for all $t\geq 0$, we have
\begin{equation}\label{21}
\mathcal{L}_\lambda(f; x, t)\longrightarrow \exp\Big(-\Exp[\bfp^{\alpha}]t^{\alpha}C(d,\beta)\Big),
\end{equation}
where $C(d, \beta)=\frac{\omega(d)}{\beta}\int_0^\infty (1-e^{-s})s^{-1-\alpha}ds$
and $\omega(d)$ is the surface area of the unit ball in $\R^d$.
\end{lemma}
Note that $C(d, \beta)$ is finite because $(1-e^{-s})\leq s$ for all $s\geq 0$.

\noi{\bf Proof:} Because of stationarity it is enough to prove the result for $x=0$. Also there is nothing to prove if $t=0$. Assume $t>0$.
First we observe from \eqref{gen-fn} that 
\begin{align}\label{105}
&\mathcal{L}_\lambda(f; 0, t)\nonumber
\\
&=\exp\Big[-\int_{\R_+}\omega(d)\int_{0}^\varrho\Big(1-e^{-\frac{\lambda^{-\kappa}t\, p}{r^\beta}}\Big)\lambda r^{d-1}\; dr F(dp)
-\int_{\R_+}\int_{\|y\|\geq \varrho}\Big(1-e^{-{\lambda^{-\kappa}tp}{f(\|y\|)}}\Big)\lambda\, m(dy) F(dp)\Big].
\end{align}
Now for the first term on the r.h.s. we use the change of variable 
$pt r^{-\beta} \lambda^{-\kappa}=s$ to obtain 
\begin{align*}
\int_{0}^\varrho\Big(1-e^{-\frac{n^{-\kappa}t\, p}{r^\beta}}\Big)\lambda r^{d-1}\; dr
= \frac{\omega(d)}{\beta}p^\alpha\, t^\alpha\int_{\lambda^{-\kappa}\varrho^{-\beta}tp}^\infty
(1-e^{-s})s^{-1-\alpha}ds.
\end{align*}
Thus by dominating convergence theorem we get as $\lambda\to \infty$, 
$$\int_{\R_+}\omega(d)\int_{0}^\varrho\Big(1-e^{-\frac{\lambda^{-\kappa}t\, p}{r^\beta}}\Big)\lambda r^{d-1}\; dr F(dp)\to
\Exp[\bfp^\alpha]t^\alpha C(d, \beta).$$
Thus we only need to show from \eqref{105} that the second term
inside the exponential goes to $0$ with $\lambda\to\infty$.
Since $(1-e^{-x})\leq x,$ for $x\geq 0,$ we have
\begin{align*}
\int_{\R_+}\int_{\|y\|\geq \varrho}\Big(1-e^{-{\lambda^{-\kappa}tp}{f(\|y\|)}}\Big)\lambda\, m(dy) F(dp)
&\leq\; \int_{\R_+}\int_{\|y\|\geq \varrho}{\lambda^{-\kappa}tp}{f(\|y\|)} 
\lambda\, m(dy) F(dp)
\\
&\leq\; t\, \lambda^{1-\kappa}\int pF(dp) \Big[\omega(d) \int_{\varrho}^\infty rf(r)dz\Big]
\\
&\to\; 0,
\end{align*}
as $\lambda\to\infty$ where we use that $\kappa>1$. Hence the proof.
\hfill $\Box$

\begin{remark}
It is interesting to note that Lemma~\ref{lem1} holds without the assumption that
$\sup_{r\geq\varrho}f(r)<\infty$.
\end{remark}

\begin{remark}
A non-degenerate random variable $X$ is said to be stable if for any $a, b>0$
and two independent copies $(X_1, X_2)$ of $X$, $aX_1+b X_2$ has the same distribution
as $cX+d$ for some $c>0$ and $d\in\R$ \cite{Em-Kl-Mi}. Stable random variables are characterized by their
characteristic functions and Laplace transforms. $X$ is said to be an $\alpha$-stable random
variable ($0<\alpha<1$) if its Laplace transform is given by $\Exp[\exp(-tX)]=\exp(-\eta|t|^\alpha)$
for $t\geq 0$ and some constant $\eta$. 
\end{remark}

The fact that stable laws show up in the context of (\ref{res-fn}) is not
new (see e.g. \cite{Baccelli-Bartek,haenggi}). What is new,
to the best of our knowledge, is 
the fact that the distribution of the finite dimensional
marginals of the limiting random field is of product form:

\begin{theorem}\label{theo-1}
Let $x_i, i=1,2,\ldots, k$ be $k$-distinct points in $\R^d$. Also let $t_i, i=1,2,\ldots, k$ be positive real numbers. By 
$\mathcal{L}_\lambda(f; x_1,\ldots, x_k, t_1,\ldots, t_k)$ we denote the Laplace transform of $\big(\hat{\I}_\lambda(f; x_1),
\ldots, \hat{\I}_\lambda(f; x_k)\big)$ at the point $(t_1, \ldots, t_k)$. Then
$$\lim_{\lambda\to\infty}\mathcal{L}_\lambda(f; x_1,\ldots, x_k, t_1,\ldots, t_k)\;=\;\exp\Big(-\Exp[\bfp^{\alpha}]C(d, \beta)\sum_{i=1}^k t^{\alpha}_i\Big), \quad \alpha=\frac{d}{\beta}. $$
\end{theorem}

In view of Theorem \ref{theo-1} we see that $\big(\hat{\I}_\lambda(f; x_1), \ldots, \hat{\I}_\lambda(f; x_k)\big)$ converges in distribution to $(\xi_1, \ldots, \xi_k)$ as $\lambda\to\infty$ where $(\xi_1, \ldots, \xi_k)$ is an i.i.d. sequence with 
$$\Exp[e^{-t \xi_1}]\;=\;\exp\big(-C(d, \beta)\Exp[\bfp^{\alpha}]t^{\alpha}\big),$$
i.e., $\xi_1$ has a $\alpha$-stable distribution. Therefore the limiting random field of $\hat{\I}_\lambda(f;\cdot)$ is a stationary random field and
its finite dimensional distributions are given by independent and identically distributed $\alpha$-stable variables.   
It makes sense to call the limiting field a
\textit{$\alpha$-stable white noise field}.

\spa

\noi{\bf Proof of Theorem \ref{theo-1}:}
Recall that $L(x, y, p)=p\, f(\|x-y\|)$. Using Proposition \ref{laplace} we have
\begin{align*}
\mathcal{L}_\lambda(x_1,\ldots, x_k; t_1,\ldots, t_k)\;=\;
\exp\Big\{-\lambda\int_{\R_+}\int_{\R^d}\Big(1-e^{-\lambda^{-\kappa}\sum_{i=1}^k t_i L(x_i, y, p)}\Big)m(dy)F(dp)\Big\}.
\end{align*}
Now for each $p$ we write
$$\Big(1-e^{-\lambda^{-\kappa}\sum_{i=1}^k t_i L(x_i, y, p)}\Big)=\sum_{i=1}^k 
e^{-\lambda^{-\kappa}\sum_{l=1}^{i-1}t_l L(x_l, y, p)}\Big(1-e^{-\lambda^{-\kappa}t_i L(x_i, y, p)}\Big).$$
Hence
\begin{align}\label{002}
&\quad \lambda\int_{\R_+}\int_{\R^d}\Big(1-e^{-\lambda^{-\kappa}\sum_{i=1}^k t_iL(x_i, y, p)}\Big)m(dy)F(dp)\nonumber
\\
&=\;\sum_{i=1}^k \lambda\int_{\R_+}\int_{\R^d} e^{-\lambda^{-\kappa}\sum_{l=1}^{i-1}t_l L(x_l, y, p)}
\Big(1-e^{-\lambda^{-\kappa}t_i L(x_i, y, p)}\Big)m(dy) F(dp).
\end{align}
%Now we calculate the $k$ terms on the r.h.s. of \eqref{002} for large $n$.
Note that the first term of the summation has already been calculated in Lemma \ref{lem1}. So we care for $i$-th terms for $i\geq 2$. Fix $i\geq 2$. Then

\begin{align}\label{003}
&\quad \lambda\int_{\R_+}\int_{\R^d} e^{-\lambda^{-\kappa}\sum_{l=1}^{i-1}t_l L(x_l, y, p)}
\Big(1-e^{-\lambda^{-\kappa}t_i L(x_i, y, p)}\Big)
m(dy) F(dp)\nonumber
\\
&\leq\; 
\lambda\int_{\R_+}\int_{\R^d} 
\Big(1-e^{-\lambda^{-\kappa}t_i L(x_i, y, p)}\Big)m(dy) F(dp)\nonumber
\\
&\to \; C(d, \beta)\Exp[\bfp^{\alpha}]t_i^{\alpha},
\end{align}
where the last line follows from Lemma~\ref{lem1}
(see for example \eqref{105}). Let $\delta=\frac{1}{2} \min_{i\neq j}\|x_i-x_j\|$. Define
$B_i=\cup_{l=1}^{i-1} B_\delta(x_l)$ where $B_\delta(x_l)$ denotes the $\delta$-open ball around $x_l$ in $\R^d$. Thus $B_i$ is the union of $\delta$-neighborhoods of $x_1,\ldots, x_{i-1}$. By definition $x_i\notin B_i$. Then for any $\eps>0$ and $M>0$ we get

\begin{align}\label{04}
& \lambda\int_{\R_+}\int_{\R^d} e^{-\lambda^{-\kappa}\sum_{l=1}^{i-1}t_l L(x_l, y, p)}
\Big(1-e^{-\lambda^{-\kappa}t_i L(x_i, y, p)}\Big)
m(dy) F(dp)\nonumber
\\
&\geq\,  \lambda\int_{\R_+}\int_{B_i^c} e^{-\lambda^{-\kappa}\sum_{l=1}^{i-1}t_l L(x_l, y, p)}
\Big(1-e^{-\lambda^{-\kappa}t_i L(x_i, y, p)}\Big)
m(dy) F(dp)\nonumber
\\
&\geq\,  (1-\eps)\lambda\int_{\{p\leq M\}}\int_{B_i^c}
\Big(1-e^{-\lambda^{-\kappa}t_i L(x_i, y, p)}\Big)
m(dy) F(dp),
\end{align}
for all large $\lambda$ where in the last line we use the fact that $L(x_l, y, p)\leq {M}\sup_{z\geq \delta} f(z)$ for $(y, p)\in B_i^c\times\{p\leq M\}, l\in\{1, \ldots, i-1\},$ and $\lambda^{-\kappa}\to 0$ as $\lambda\to\infty$.

Again on $B_i$, we have $L(x_i, y, p)\leq p \sup_{z\geq \delta} f(z)$. Hence using the fact
$(1-e^{-x})\leq x$ for all $x\geq 0$, we get
\begin{align*}
\lambda\int_{\{p\leq M\}}\int_{B_{i}}\Big(1-e^{-\lambda^{-\kappa}t_iL(x_i, y, p)}\Big)
m(dy) F(dp)\leq 
\lambda\, [\sup_{z\geq \delta} f(z)]\, m(B_\delta(0))\int_{\{p\leq M\}}t_i\, p \lambda^{-\kappa}F(dp).
\end{align*}
Since $\kappa=\frac{\beta}{d}>1$ we see that the r.h.s. of the above expression tends to $0$ as $\lambda\to\infty$ for every fixed $M$.
Thus combining with \eqref{04} we see that
\begin{align*}
&\liminf_{\lambda\to\infty}\;  \lambda\;\int_{\R_+}\int_{\R^d}
 e^{-\lambda^{-\kappa}\sum_{l=1}^{i-1}t_l L(x_l, y, p)}
\Big(1-e^{-\lambda^{-\kappa}t_i L(x_i, y, p)}\Big)
m(dy) F(dp)
\\
&\geq\; \liminf_{\lambda\to\infty}\quad (1-\eps)
\lambda\int_{\{p\leq M\}}\int_{\R^d}\Big(1-e^{-\lambda^{-\kappa}t_i L(x_i, y, p)}\Big)
m(dy) F(dp),
\\
&=\;(1-\eps)\int_{\{p\leq M\}}p^{d/\beta}C(d, \beta)t_i^{d/\beta} F(dp),
\end{align*}
for any $\eps, M>0$ where for the last equality we can follow the same computation as in Lemma \ref{lem1}.
Now using the fact that $\Exp[\bfp^{d/\beta}]<\infty$, we let $M\to\infty$ and $\eps\to 0$ to obtain
\begin{align}\label{004}
\liminf_{\lambda\to\infty}\; \lambda\; \int_{\R_+}\int_{\R^d} 
e^{-\lambda^{-\kappa}\sum_{l=1}^{i-1}t_l L(x_l, y, p)}
\Big(1-e^{-\lambda^{-\kappa}t_i L(x_i, y, p)}\Big)
m(dy) F(dp)
\;\geq\; \Exp[\bfp^{d/\beta}]C(d, \beta)t_i^{d/\beta}.
\end{align}
Hence combining \eqref{003} and \eqref{004} we get for $i\geq 2$,
\begin{align*}
\lim_{\lambda\to\infty}\; \lambda\;\int_{\R_+}\int_{\R^d} 
e^{-\lambda^{-\kappa}\sum_{l=1}^{i-1}t_l L(x_l, y, p)}
\Big(1-e^{-\lambda^{-\kappa}t_i L(x_i, y, p)}\Big)
dy F(dp)
\;=\; \Exp[\bfp^{\alpha}]C(d, \beta)t_i^{\alpha}.
\end{align*}
Thus the result follows from \eqref{002}.\hfill $\Box$

\subsection{Scaling Limits of Extremal Shot-noise Fields}

This subsection is devoted to the analysis
of the \textit{extremal-random field} generated by the
response function \eqref{res-fn} (and also by \eqref{gen-fn},
see Remark~\ref{rem3.2}).
Let $\tilde\Phi$ be given independently marked PPP on $\R^d\times (0, \infty)$. We define the \textit{extremal-random field} at a
point $y$ as follows:
\begin{equation}\label{max-fl}
\calM(y)=\sup_{(x_i, \bfp_i)\in\tilde\Phi}\frac{\bfp_i}{\|x_i-y\|^\beta}, \quad \beta>d.
\end{equation}
We see that $\calM(y)\leq \I(y)$. Therefore $\calM$ is an almost surely finite and stationary random field. When the underlying PPP
has intensity $\lambda$ we denote the extremal-random field by $\calM_\lambda$. Define
$$\hat{\calM}_\lambda=\frac{1}{\lambda^\kappa}\calM,\quad \text{for}\, \, 
\kappa=\frac{\beta}{d}.$$
Let $x_1, x_2, \ldots, x_k$ be $k$-distinct points in $\R^d$. By $\calG_\lambda$ we denote the multivariate cumulative distribution of
$(\hat{\calM}_\lambda(x_1), \ldots, \hat{\calM}_\lambda(x_k))$ , i.e.,
$$\calG_\lambda(t_1, \ldots, t_k)=\Pb\Big[\hat{\calM}_\lambda(x_1)\leq t_1, \ldots, \hat{\calM}_\lambda(x_k)\leq t_k\Big],$$
for $t_i\geq 0,\, i=1,\ldots, k$.
The proof of the following lemma is similar to that
of \cite[Proposition~2.4]{dombry}.
\begin{lemma}\label{lem-2}
Let $(t_1,\ldots, t_k)\in [0, \infty)^k$. Then
$$\calG_\lambda(t_1, \ldots, t_k)=\exp\Big\{-\lambda\int_{\R^d}\int_0^\infty\Big(1-\Pi_{i=1}^k 1_{\{\frac{p}{\|y-x_i\|^\beta}\leq \lambda^{\kappa}t_i\}}\Big)
F(dp)\, m(dy)\Big\},$$
where $\calG_\lambda$ is defined above.
\end{lemma}

\noindent{\bf Proof:} For simplicity, we prove the lemma for $\lambda=1$. 
The proof for general $\lambda$ is analogous. First we observe that
$$1_{\{\calM(x_i)\leq t_i\}} =\Pi_{(y_j, \bfp_j)\in\tilde\Phi}1_{\{\frac{\bfp_j}{\|y_j-x_i\|^\beta}\leq t_i\}}.$$
Therefore 
\begin{align*}
\Pb\Big[\calM(x_1)\leq t_1, \ldots, \calM(x_k)\leq t_k\Big]
&= \Exp\Big[\Pi_{(y_j, \bfp_j)\in\tilde\Phi}1_{\{\frac{\bfp_j}{\|y_j-x_1\|^\beta}\leq t_1\}}\cdots
\Pi_{(y_j, \bfp_j)\in\tilde\Phi}1_{\{\frac{\bfp_j}{\|y_j-x_k\|^\beta}\leq t_k\}}\Big]
\\
& =\; \Exp[e^{H}],
\end{align*}
where 
$$H=\sum_{(y_j, \bfp_j)\in\tilde\Phi}\log\big(1_{\{\frac{\bfp_j}{\|y_j-x_1\|^\beta}\leq t_1\}}\big)+\cdots+
\sum_{(y_j, \bfp_j)\in\tilde\Phi}\log\big(1_{\{\frac{\bfp_j}{\|y_j-x_k\|^\beta}\leq t_k\}}\big).$$
Therefore from Proposition~\ref{laplace} we have
\begin{align*}
\calG_1(t_1, \ldots, t_k)& =\;\exp\Big\{-\int_{\R^d}\bigg(1-\int_0^\infty \exp\Big[\sum_{i=1}^k\log(1_{\{\frac{p}{\|y-x_i\|^\beta}\leq t_i\}})\Big]F(dp)\bigg)m(dy)\Big\}
\\
&=\; \exp\Big\{-\int_{\R^d}\int_0^\infty \Big(1-\Pi_{i=1}^k1_{\{\frac{p}{\|y-x_i\|^\beta}\leq t_i\}} \Big)F(dp)m(dy)\Big\}.
\end{align*}
Hence the proof. \hfill $\Box$

\begin{theorem}\label{theo-ex}
Let $x_1, \ldots, x_k$ be $k$-distinct points in $\R^d$. Consider $\calG_\lambda$ as defined above.
Let 
$$\gamma(d,\beta):=\frac{\omega(d)}{\beta}\int_0^\infty \Pb(p>s)s^{-1+d/\beta}ds=\frac{\omega(d)}{d} \Exp[p^{d/\beta}]<\infty,$$
where $\omega(d)$ denotes the surface area of the unit ball in $\R^d$.
 Then
$$\lim_{\lambda\to\infty}\calG_\lambda(t_1, \ldots, t_k)\;=\;\Pi_{i=1}^k \exp\Big(-\gamma(d, \beta)t_i^{-\alpha}\Big),\quad \alpha=\frac{d}{\beta},$$
for $t_i\geq 0, i=1, \ldots, k$.
\end{theorem}

By Theorem~\ref{theo-ex} we see that the limiting random field of the extremal field is a max stable random field, i.e., the
finite dimensional distribution of the limiting fields are given by collection of i.i.d. \textit{Fr\'{e}chet} distribution
of exponent $\alpha=d/\beta$ (\cite{Em-Kl-Mi}). We call the limiting field a \textit{$\alpha$-Fr\'{e}chet white noise field}.  
This result is similar to that obtained in \cite[Section~3.2]{dombry}. The results in \cite{dombry} are obtained for
$f$ that are integrable on $\mathbb{R}^d$ and the finite dimensional distributions of limiting random field are not necessarily
independent there. 

\noindent{\bf Proof:} The proof is similar to the proof of Theorem~\ref{theo-1}.
However, we add it for clarity.
Without loss of generality we assume that $t_i>0$ for all $i=1, \ldots, k$. Now by Lemma~\ref{lem-2} we have
\begin{equation}\label{31}
\calG_\lambda(t_1, \ldots, t_k)\;=\;\exp\Big\{-\lambda\int_{\R^d}\int_0^\infty\Big(1-\Pi_{i=1}^k 1_{\{\frac{p}{\|y-x_i\|^\beta}\;\leq\; \lambda^{\kappa}t_i\}}\Big)
F(dp)\, m(dy)\Big\}.
\end{equation}
As earlier (see the display preceding \eqref{002}) we write 
\begin{align*}
\Big(1-\Pi_{i=1}^k 1_{\{\frac{p}{\|y-x_i\|^\beta}\;\leq\; \lambda^{\kappa}t_i\}}\Big) 
&=\;\sum_{i=1}^k \Pi_{l=1}^{i-1}
1_{\{\frac{p}{\|y-x_l\|^\beta}\;\leq\; \lambda^{\kappa}t_l\}}
\Big(1-1_{\{\frac{p}{\|y-x_i\|^\beta}\;\leq\; \lambda^{\kappa}t_i\}}\Big)
\\
& =\; \sum_{i=1}^k \Pi_{l=1}^{i-1}1_{\{\frac{p}{\|y-x_l\|^\beta}\;\leq\; 
\lambda^{\kappa}t_l\}}
\Big(1_{\{\frac{p}{\|y-x_i\|^\beta}\;>\; \lambda^{\kappa}t_i\}}\Big).
\end{align*}
We notice that to prove the theorem we only need to prove the convergence of the exponent in \eqref{31}. Thus we write
\begin{align}\label{32}
&\quad \lambda\int_{\R^d}\int_0^\infty\Big(1-\Pi_{i=1}^k 1_{\{\frac{p}{\|y-x_i\|^\beta}\;\leq\; \lambda^{\kappa}t_i\}}\Big)
F(dp)\, m(dy) \nonumber
\\
&=\;\sum_{i=1}^k \lambda\int_{\R^d}\int_0^\infty\Pi_{l=1}^{i-1}1_{\{\frac{p}{\|y-x_l\|^\beta}
\;\leq\; \lambda^{\kappa}t_l\}}
\Big(1_{\{\frac{p}{\|y-x_i\|^\beta}\;>\; \lambda^{\kappa}t_i\}}\Big)F(dp)\, m(dy).
\end{align}
A simple change of variable in \eqref{32} shows that 
\begin{equation}\label{33}
\lambda\int_{\R^d}\int_0^\infty\Big(1-\Pi_{i=1}^k 1_{\{\frac{p}{\|y-x_i\|^\beta}\;\leq\; \lambda^{\kappa}t_i\}}\Big)
F(dp)\, m(dy)
\;\leq\; \sum_{i=1}^k\gamma(d, \beta) t_i^{-d/\beta}.
\end{equation}
Now recall the sets $B_i$ from the proof of Theorem~\ref{theo-1}. Then for any positive $M$ we have, for every $i$,
\begin{align*}
&\liminf_{\lambda\to\infty}\; \lambda\;\int_{\R^d}\int_0^\infty \Pi_{l=1}^{i-1}1_{\{\frac{p}{\|y-x_l\|^\beta}\;\leq\; \lambda^{\kappa}t_l\}}
\Big(1_{\{\frac{p}{\|y-x_i\|^\beta}\;>\; \lambda^{\kappa}t_i\}}\Big)F(dp)\, m(dy)
\\
&\geq\; \liminf_{\lambda\to\infty}\; \lambda\; \int_{B_i^c}\int_0^M \Pi_{l=1}^{i-1}1_{\{\frac{p}{\|y-x_l\|^\beta}\;\leq\; \lambda^{\kappa}t_l\}}
\Big(1_{\{\frac{p}{\|y-x_i\|^\beta}\;>\; \lambda^{\kappa}t_i\}}\Big)F(dp)\, m(dy)
\\
&=\; \liminf_{\lambda\to\infty}\; \lambda\; \int_{B_i^c}\int_0^M 
\Big(1_{\{\frac{p}{\|y-x_i\|^\beta}\;>\; \lambda^{\kappa}t_i\}}\Big)F(dp)\, m(dy),
\end{align*}
where we used the fact that  $\sup_{l\leq i-1}\sup_{y\in B_i^c}\frac{1}{\|y-x_l\|^\beta}\;<\infty$. Since 
$\text{dist}(x_i, B_i)\geq\frac{1}{2}\min_{i\neq j}\|x_i-x_j\|$ we get
\begin{align*}
&\quad \liminf_{\lambda\to\infty} \;\lambda\;\int_{\R^d}\int_0^\infty \Pi_{l=1}^{i-1}1_{\{\frac{p}{\|y-x_l\|^\beta}\;\leq\; \lambda^{\kappa}t_l\}}
\Big(1_{\{\frac{p}{\|y-x_i\|^\beta}\;>\; \lambda^{\kappa}t_i\}}\Big)F(dp)\, m(dy)
\\
&\geq\;  \liminf_{\lambda\to\infty} \;\lambda \;\int_{\R^d}\int_0^M 
\Big(1_{\{\frac{p}{\|y-x_i\|^\beta}\;>\; \lambda^{\kappa}t_i\}}\Big)F(dp)\, m(dy),
\\
&=\;\frac{\omega(d)}{\beta}t^{-d/\beta}_i\int_0^\infty \Pb\big(s<\bfp\leq M\big)s^{-1+d/\beta}ds.
\end{align*}
Now let $M\to\infty$ to obtain
\begin{align}\label{34}
&\quad \liminf_{\lambda\to\infty} \lambda\int_{\R^d}\int_0^\infty \Pi_{l=1}^{i-1}1_{\{\frac{p}{\|y-x_l\|^\beta}\;\leq\; \lambda^{\kappa}t_l\}}
\Big(1_{\{\frac{p}{\|y-x_i\|^\beta}\;>\; \lambda^{\kappa}t_i\}}\Big)F(dp)\, m(dy)\nonumber
\\
&\geq\; \frac{\omega(d)}{\beta}t^{-d/\beta}_i\int_0^\infty \Pb(\bfp\;>\;s)s^{-1+d/\beta}ds=\gamma(d,\beta)t^{-d/\beta}_i.
\end{align}
The proof is completed by combining \eqref{31}, \eqref{33} and \eqref{34}. \hfill $\Box$

\begin{remark}\label{rem3.2}
The result of Theorem~\ref{theo-ex} extend to
extremal random fields defined using any positive $f$
satisfying the conditions in \eqref{gen-fn}.
\end{remark}

\section{Applications to Stochastic Wireless Networks}
\label{secapp}
\subsection{Scaling of the Shannon Rate with Distance}
\label{ssec:distan}
In this section we give a first application of the above results to
SIR stochastic models \cite{Baccelli-Bartek, haenggi}.
In the rest of this article we will consider $d=2$. We start by defining
the Signal to Interference Ratio (SIR), which finds its root in Shannon's
Channel Coding Theorem \cite{CoT91}. Let $\Phi$ be a given PPP.
The support of $\Phi$ represents the network nodes on the plane. We consider
two fixed additional points in $\R^2$: $0$ and $e_\lambda=(\lambda, 0)$.
Let $\{F_{0\lambda}, F_{i\lambda} ; i\geq 1\}$ be a collection of positive i.i.d. random variables.
The variable $F_{i\lambda}$ represents a random perturbation
called the \textit{fading} from $x_i\in\Phi$ to $e_\lambda$.
$F_{0\lambda}$ denotes the fading between $0$ and $e_\lambda$.
We assume that the fading is independent of the PPP.
Let $F(dp)$ be the common distribution of the fading variables.
Let $\ell:[0, \infty)\to[0, \infty]$ be the \textit{path-loss function}
given by
$$\ell(r)=r^\beta, \quad \text{for some}\ \beta>2.$$ 
We also assume that
\begin{equation}\label{10}
\int_{\R_+}p\, F(dp)<\infty.
\end{equation}
Node $e_\lambda$ receives the signal from $x_i\in\Phi$
with power $\frac{F_{i\lambda}}{\ell(\|e_\lambda-x_i\|)}$.
Hence the total power received at $e_\lambda$
from $\Phi$ is given by the shot-noise
\begin{equation}\label{11}
\I(e_\lambda)=\sum_{x_i\in\Phi} L(e_\lambda, x_i, F_{i\lambda}), \quad \text{with}\quad L(e_\lambda, x, p)=\frac{p}{\ell(\|e_\lambda-x\|)}.
\end{equation}

Now we define the SIR between $0$ and $e_\lambda$ as
$$\SINR_{0\lambda}=\cfrac{F_{0\lambda}/\ell(\|e_\lambda\|)}{\I(e_\lambda)}.$$
It follows from Campbell's theorem 
\cite{SKM} that $\I(e_\lambda)$ is finite with probability 1 so that
$\SINR_{0\lambda}$ is positive with probability $1$.
Let $c>0$.
From Shannon's Channel Coding Theorem,
when treating interference as Gaussian noise,
the transmission from 0 to $e_\lambda$ is possible at
rate $\frac{1}{2} \ln (1+c)$ if
$\SINR_{0\lambda}>c$ and impossible if $\SINR_{0\lambda}<c$.
This Shannon rate explains the practical
importance of the following scaling result:

% if $e_n\notin \text{support}(\Phi)$ with probability $1$. This is indeed the case when $\Phi$ is a stationary PPP.

\begin{theorem}\label{theo-2}
Let $\Phi$ be a stationary PPP of intensity $1$.
 Then under Condition \eqref{10} we have
\begin{equation}
\liminf_{c\to 0+}\;\liminf_{\lambda\to\infty}
\;\Pb\Big(\SINR_{0\lambda}\geq \frac{c}{\lambda^{\beta}}\Big)
\;=\;1.
\end{equation}
\end{theorem}

\noi{\bf Proof:} The proof is based on the following observation:
if $\Phi_1$ is a stationary PPP of intensity $1$, then the PPP 
obtained by the mapping $x\mapsto \frac{x}{\lambda}$ has
intensity $\lambda^2$.
Let $\Phi_{\lambda^2}$ be the stationary PPP of intensity $\lambda^2$.
Therefore we see from \eqref{11} that $\I(e_n)$ has the same distribution
as $\frac{1}{\lambda^\beta}\I_{\lambda^2}(e_1)$, where
$$\I_{\lambda^2}(e_1)=\sum_{x_i\in\Phi_{\lambda^2}} L(e_1, x_i, F_{i\lambda}).$$
Let $\xi$ be an $\alpha$-stable random variable such that for any $t\geq 0$,
$$ \Exp[e^{-t \xi}]=\exp(-t^{\alpha} C(\beta)),$$
where $C(\beta)=C(2, \beta)\int_{\R_+}p^{\alpha}\, F(dp), \alpha=2/\beta,$ and $C(\beta)$ is given by Lemma \ref{lem1}.
 By Lemma \ref{lem1} we know that
$ \frac{1}{\lambda^\beta}\I_{\lambda^2}(e_1)\to\xi$ in the sense of
convergence in distribution. Therefore for any $c>0$,
we obtain
\begin{align*}
\Pb\Big(\SINR_{0\lambda}\geq \frac{c}{\lambda^{\beta}}\Big)&=\;\Pb\Big(\frac{F_{0\lambda}/\ell(\|e_1\|)}{\I_{\lambda^2}(e_1)}\geq 
\frac{c}{\lambda^{\beta}}\Big)
\\
&=\;\Pb\Big(\frac{F_{0\lambda}}{c\ell(\|e_1\|)}\geq \lambda^{-\beta}\I_{n^2}(e_1)\Big)
\\
&\geq\; \Pb\Big(\frac{F_{0\lambda}}{c\ell(\|e_1\|)}\geq \lambda^{-\beta}\I_{\lambda^2}(e_1), F_{0\lambda}>\delta\Big)
\\
&\geq\; \Pb\Big(\frac{\delta}{c\ell(\|e_1\|)}\geq \lambda^{-\beta}\I_{\lambda^2}(e_1), F_{0\lambda}>\delta\Big)
\end{align*}
for any positive constant $\delta$. Therefore using the independence we have
\begin{align*}
\Pb\Big(\SINR_{0\lambda}\geq \frac{c}{\lambda^{\beta}}\Big) &\geq\; \Pb\Big(\frac{\delta}{c\ell(\|e_1\|)}\geq \lambda^{-\beta}\I_{\lambda^2}(e_1)\big)\Pb
\big( F_{0\lambda}>\delta\Big)
\\
& =\;\Pb\Big(\frac{\delta}{c\ell(\|e_1\|)}\geq \lambda^{-\beta}\I_{\lambda^2}(e_1)\Big)\int_{\R_+}1_{\{p>\delta\}}F(dp).
\end{align*}
Therefore letting $\lambda\to\infty$, we get
$$\liminf_{\lambda\to\infty}\; \Pb\Big(\SINR_{0\lambda}\geq\;
 \frac{c}{\lambda^{\beta}}\Big)\;\geq\;
 \Pb\Big(\frac{\delta}{c\ell(\|e_1\|)}> \xi\Big)\int_{\R_+}1_{\{p>\delta\}}F(dp),$$
where $\xi$ is an $\alpha$-stable random variable.

Now first let $c\to 0+$ and then $\delta\to 0$ to obtain the result. \hfill $\Box$

Another way of rephrasing Theorem \ref{theo-2} is that 
for a Poisson field of interferers with density $\lambda$,
the scale at which the $\SINR$ decreases for
a link of length 1 is $\lambda^{-\kappa}, \, \kappa=\frac{\beta}{2},$ and the Shannon rate
on that link scales like $\lambda^{-\kappa}$.

\subsection{SINR Percolation in Ultra Dense Networks}
Signal to Interference and Noise Ratio
(SINR) Percolation received a lot of attention
(see \cite{Fran-Mees,Baccelli-Bartek} and the references therein)
but was only studied in the case of a bounded 
response function to the best of our knowledge.

The aim of this subsection is to discuss SINR Percolation
for the power law response functions considered here
in terms of scaling laws for ultra-dense networks.

Let $v_1, \ldots, v_k$ be $k$-given distinct points on $\R^2$.
Here we are interested in finding the
scale of SINR at these points when the network density tends to infinity.
Let $\Phi=\Phi_\lambda$ be a stationary PPP of intensity $\lambda$.
Let $(W_1, \ldots, W_k)$ be $k$ non-negative i.i.d. random
variables that are independent of $\Phi_\lambda$. We may think of
$W_i$ as the power of \textit{thermal noise} at $v_i$.
Let $\{F_{l(l+1)}, F_{ij}; i, j, n\geq 1\}$ be a family
of non-negative distributions with common cumulative
distribution $F(dp)$ where $F(dp)$ satisfies Condition \eqref{10}.
As earlier we may think of
$F_{l({l+1})}$ as the fading variables
between nodes at $v_l$ and $v_{l+1}$. Let
$F_{ij}$ denote the fading variable between $x_i\in\Phi_\lambda$ and $v_j$.
As earlier we assume that the fading is independent of
$\Phi_\lambda$ and $\{W_i\}$. We define the SINR 
between $v_l$ and $v_{l+1}$ as follows
$$ S_{l(l+1)}(\lambda)=\SINNR_{l(l+1)}=\frac{F_{l(l+1)}/\ell(\|v_l-v_{l+1}\|)}{W_{l+1}+\I_\lambda(v_{l+1})},$$
where
$$\I_\lambda(v_{l+1})=\sum_{x_i\in\Phi_\lambda} L(v_{l+1}, x_i, F_{i(l+1)}) \quad \text{and}\quad L(v, x, p)=\frac{p}{\ell(\|v-x\|)}.$$

\begin{theorem}\label{theo-3}
Let $c>0$ be given. Let $\gamma=\min_{1\leq i\leq k-1}\frac{1}{\ell(\|v_i-v_{i+1}\|)}$. Then 
\begin{align*}
&\liminf_{\lambda\to\infty}\Pb(S_{12}(\lambda)\geq c\lambda^{-\kappa},\ldots,
 S_{(k-1)k}(\lambda)\geq c\lambda^{-\kappa})
\\
&\geq 
\Big[\Pb(\xi<\frac{\gamma}{2\sqrt{c}})\big(\int_{\{p\geq\sqrt{c}\}}F(dp)\big)\Pb(W_2\leq\frac{1}{\sqrt{c}})\Big]^{k-1},
\end{align*}
where $\xi$ is a $\alpha$-stable variable such that $\Exp[e^{-t\xi}]=\exp(-C(\beta) t^{\alpha})$ for $t\geq 0$ and $C(\beta)=C(2, \beta)$ is given by
Theorem \ref{theo-2}. In particular, 
$$\lim_{c\to 0+}\liminf_{\lambda\to\infty}\Pb(S_{12}(\lambda)\geq c \lambda^{-\kappa},\ldots, S_{(k-1)k}(\lambda)\geq c\lambda^{-\kappa})=1.$$
\end{theorem}

\noi{\bf Proof:} Let $\gamma=\min_{1\leq i\leq k-1}\frac{1}{\ell(\|v_i-v_{i+1}\|)}$. Then for any $c>0$
\begin{align*}
& \Pb\Big(S_{12}(\lambda)\geq c\lambda^{-\alpha},\ldots, S_{(k-1)k}(\lambda)
\geq c\lambda^{-\alpha}\Big)
\\
& \geq \Pb\Big(\lambda^{-\alpha}\I_\lambda(v_2)\leq \big(\frac{\gamma}{c}F_{12}-\lambda^{-\alpha}W_2\big), \ldots,
\lambda^{-\alpha}\I_\lambda(v_k)\leq \big(\frac{\gamma}{c}F_{(k-1)k}
-\lambda^{-\alpha}W_k\big)\Big)
\\
&\geq \Pb\Big(\lambda^{-\alpha}\I_\lambda(v_2)\leq \big(\frac{\gamma}{c}F_{12}-\lambda^{-\alpha}W_2\big), \ldots,
\lambda^{-\alpha}\I_\lambda(v_k)\leq \big(\frac{\gamma}{c}F_{(k-1)k}-
\lambda^{-\alpha}W_k\big), A(\delta)\Big)
\end{align*}
for any positive $\delta$ where $A(\delta)=\{ F_{i(i+1)}\geq \delta, \forall\ i=1,\ldots, k-1\}\cap\{W_{i}<1/\delta, \forall\ i=2, \ldots , k\}$.
Therefore 
\begin{align*}
& \Pb\big(S_{12}(\lambda)\geq c\lambda^{-\kappa},\ldots, S_{(k-1)k}(\lambda)\geq
 c\lambda^{-\kappa}\big)
\\
&\geq \Pb\Big(\lambda^{-\kappa}\I_\lambda(v_2)\leq \big(\frac{\gamma}{c}\delta
-\lambda^{-\kappa}\frac{1}{\delta}), \ldots,
\lambda^{-\kappa}\I_\lambda(v_k)\leq \big(\frac{\gamma}{c}\delta-\lambda^{-\kappa}\frac{1}{\delta}\big), A(\delta)\Big)
\\
&\geq \Pb\Big(\lambda^{-\kappa}\I_\lambda(v_2)\leq \frac{\gamma}{2c}\delta, \ldots,
\lambda^{-\kappa}\I_\lambda(v_k)\leq \frac{\gamma}{2c}\delta\Big)\Pb\Big(A(\delta)\Big),
\end{align*}
for all large $\lambda$, where we used independence in the last line.
Therefore, applying Theorem \ref{theo-1}, we get
$$\liminf_{\lambda\to\infty}\Pb\big(S_{12}(\lambda)\geq c\lambda^{-\kappa},\ldots, S_{(k-1)k}(\lambda)\geq c\lambda^{-\kappa}\big)
\geq \Pb(\xi_1<\frac{\gamma}{2c}\delta, \ldots, \xi_{k-1}<\frac{\gamma}{2c}\delta)\Pb\Big(A(\delta)\Big),$$
where $(\xi_i, \ldots, \xi_{k-1})$ is an i.i.d.
sequence of $\alpha$-stable random variables.
The proof follows by choosing $\delta=\sqrt{c}$
and using the independence property of 
$\{F_{i(i+1)}, i=1, \ldots, k-1\}$
and $\{W_i, i=2, \ldots, k\}$.\hfill $\Box$

%{\color{red}We also need to add comments about the above theorem}
This is a much stronger result which leverages the independence of the limiting fields to assess the joint scale of decrease of the SINR threshold or the joint Shannon rate on a {\em collection of links} as above.

\vspace{.2in}

Finally we use our scaling limit result to produce results on percolation. By $\Z^2$ we denote the integer lattice on $\R^2$.
The nodes of $\Z^2$ will be considered as sites. By $(z_1, z_2)$ we denote the coordinate of points in $\Z^2$. Let
$\{g_{z\bar{z}}\; :\; z\in\Z^2, \bar{z}\in\Z\}$ be a collection of positive i.i.d. random variables
where $g_{z\bar{z}}$ denotes the fading
between the sites $z$ and $\bar{z}$. Now we consider an independently marked PPP $\tilde\Phi_\lambda$ on $\R^2$ of intensity $\lambda$.
We assume that the marks are distributed according to distribution $F(dp)$ and
$$\int_0^\infty p F(dp)<\infty.$$
Let $f$ be a measurable function with compact support and satisfying the conditions in \eqref{gen-fn}. We define the shot-noise
at location $z$ by
$$\I_\lambda(f;z)=\sum_{(x_i, \bfp_i)\in\tilde\Phi_\lambda} \bfp_i\, f(\|x_i-z\|).$$
Note that $f$ is not bounded on $[0, \infty)$. The SINR between $z$ and $\bar{z}$ is defined by
$$S_\lambda(z, \bar{z}):=\SINNR_{z\bar{z}}=\frac{g_{z\bar{z}}/f(1)}{W_{\bar{z}}
+\I_\lambda(f;\bar{z})},$$
where $\{W_{z}\; : z\in\Z^2\}$ is an independent sequence
of non-negative i.i.d. random variables representing thermal noise.
We also assume that
$\{g_{z\bar{z}}\; :\; z\in\Z^2, \bar{z}\in\Z^2\}$ is independent of $\tilde\Phi_\lambda$ for all $\lambda$. Let $\tau=\tau_\lambda>0$ denote the SINR
threshold. As explained above, only link that experience
a SINR above this threshold are operational. 
We hence say the site $z$
is connected to the site $\bar{z}$ if $S_n(z, \bar{z})\geq \tau_\lambda$.
We construct a random
graph on $\Z^2$ as follows: we put an edge between $z$ and $\bar{z}$
if and only if $z, \bar{z}$ are connected to each other and
$z, \bar{z}$ are neighbors to each other.
Note that this forms a random graph of $\Z^2$.
We denote by $\mathcal{C}^\lambda_0$ the largest connected component in the random
graph containing $0$. Let also $|\mathcal{C}^\lambda_0|$ denote the number
of sites in $\mathcal{C}^\lambda_0$. In the next theorem we show that, for
large $\lambda$ and suitably chosen $\tau$, the graph percolates.

\begin{theorem}\label{theo-4}
Let $f$ be a function with compact support satisfying
the conditions in \eqref{gen-fn}. Let $\tau_\lambda=\frac{c}{\lambda^\kappa}$
for $\kappa=\beta/2$. Then there exist positive constants $\lambda_0, c_0,$
such that, for any $\lambda\geq \lambda_0, c\in (0, c_0]$,
we have $\Pb(|\mathcal{C}^\lambda_0|=\infty)>0$.
\end{theorem}

By the above theorem we see that if we choose $c$ small enough, then
in a sufficiently dense network the site $0$ can
send/receive \textit{information} to/from infinitely many sites
with positive probability. One should also compare this result
with \cite[Theorem 2.7.4]{Fran-Mees} where percolation is established
for a bounded $f$.

\noindent{\bf Proof:}\ We say a site $z\in\Z^2$ is occupied if and only if
$$\frac{\min(g_{z^{-1}z}, g_{z^{+1}z}, g_{z_{-1}z}, g_{z_{+1}z})/f(1)}
{W_{z}+\I_n(f;z)}\geq \frac{c}{\lambda^\kappa},$$
where $z^{\pm}=(z_1, z_2)^{\pm}=(z_1\pm 1, z_2)$ and $z_{\pm}=(z_1, z_2)_{\pm}=(z_1, z_2\pm 1)$. Define
$$p\;=\;p_\lambda\;=\;\Pb\Big(\frac{\min(g_{z^{-1}z}, g_{z^{+1}z}, g_{z_{-1}z}, g_{z_{+1}z})/f(1)}
{W_{z}+\I_\lambda(f;z)}\geq \frac{c}{\lambda^\kappa}\Big).$$
Therefore every site is occupied with probability $p_\lambda$. We put an edge between two neighboring sites if and only if the 
sites are occupied.
It is easy to see that the probability that a site is occupied is the same for all sites. It is also
important to see that the random graph is stationary. But two neighboring sites might not be independent due to
the dependencies through the shot-noise field. Let the support of $f$ lie in a box of size $m$. By box of size $m$ we mean collection of
all vertices with graph distance less than $m$ from $0$.
By $D_m(z)$ we denote the box of size $m$
with center at $z\in\Z^2$. Therefore we see that
$$\I_\lambda(f;z)=\sum_{(x_i, p_i)\in\tilde\Phi, x_i\in D_m(z)}p_if(|x_i-z|).$$
Therefore, if $z$ and $\bar{z}$ are such that
$D_m(z)\cap D_m(\bar{z})=\emptyset$,
we have $(\I_\lambda(f;z), \I_\lambda(f;\bar{z}))$ i.i.d. since the
underlying point process is Poisson. Therefore the states of sites
that are at distance more than $(2m+2)$ are independent. Now
from \cite[Theorem 2.3.1]{Fran-Mees} we have a constant $p=p(m)<1$
such that if $p_\lambda>p(m)$, then $\Pb(|\bar{\mathcal{C}}^\lambda_0|=\infty)>0$,
where $\bar{\mathcal{C}}^\lambda_0$ denotes the maximal connected
component in our new random graph containing $0$. We note that 
$\bar{\mathcal{C}}_0^\lambda\subset\mathcal{C}^\lambda_0$ and therefore 
$$ \Pb\Big(|\bar{\mathcal{C}}^\lambda_0|=\infty\Big)>0\;\Rightarrow\;
 \Pb\Big(|\mathcal{C}^\lambda_0|=\infty\Big)>0.$$
Thus to complete the proof we only need to find $n$ and $c$ so that that $p_\lambda>p(d)$. We assume that $f(1)>0$ otherwise there is nothing to
prove. We note that for $c>0$, (from calculation as in Theorem~\ref{theo-3})
\begin{align*}
p_\lambda &=\;\Pb\Big(\frac{\min(g_{0^{-1}0}, g_{0^{+1}0}, g_{0_{-1}0}, g_{0_{+1}0})/f(1)}
{W_{0}+\I_\lambda(f;0)}\geq \frac{c}{\lambda^\kappa}\Big)
\\
& \geq\; \Pb(\hat{\I}_\lambda(f;0)<\frac{1}{2f(1)\sqrt{c}})\Pb(\min(g_{0^{-1}0}, g_{0^{+1}0}, g_{0_{-1}0}, g_{0_{+1}0})\geq \sqrt{c})
\Pb(W_0\leq \frac{\lambda^\kappa}{2f(1)\sqrt{c}}).
\end{align*}
Therefore, using Lemma~\ref{lem1}, we can find $\lambda_0, c_0$
such that, for $\lambda\geq \lambda_0, c\in (0, c_0]$, we have $p_\lambda>p(m)$.\hfill $\Box$

\subsection{Scaling of the Delay with Distance and an Open Question}

Another way of looking at the result of Subsection \ref{ssec:distan}
is in terms of \emph{delay}. For an $o(1)$ SIR, the 
number of bits transmitted in an $o(1)$ interval is $o(1)$.
We will then say that an $o(1)$ SIR has an $o(1)$ delay.
Theorem~\ref{theo-2} tells us that for a link of length $\lambda$
in a Poisson field of interferers
with density 1, the Shannon rate tends to 0 like $\lambda^{-\beta}$.
Hence, this theorem can be rephrased by saying
that the delay $D_\lambda$ to transmit $o(1)$ bits over distance $\lambda$ 
in one hop (in the above scheme, one sends these $o(1)$
bits directly from 0 to $\lambda$ in one hop)
scales like $\lambda^{\beta}$ as $\lambda\to\infty$.
The general case (including the possibility of multiple hops)
was studied in \cite{Baccelli-Bartek-Omid}
where it was shown that for all possible schemes within
this framework (single or multi-hop),
$\frac{\Exp[D_\lambda]}{\lambda}\to\infty$ as $\lambda\to\infty$.
This raises the following question: does there exists a
scheme that allows one to transmit $o(1)$ bits
with a delay $D_\lambda$ such that $\frac{D_\lambda}{\lambda^\gamma}$
has a non-degenerate limit for some $\gamma\in(1, \beta)$?

\vspace{.2in}

\noi\textbf{Acknowledgment:} Authors are grateful to the referees
for their fruitful comments and pointing out some important references from literature.

\end{document}